\documentclass[a4paper,10pt]{amsart}
\usepackage{amsmath,amsthm,amssymb,latexsym,enumerate,color,hyperref}
\usepackage{graphicx}

\begin{document}

\newtheorem{thm}{Theorem}
\newtheorem{prop}[thm]{Proposition}
\newtheorem{lem}[thm]{Lemma}
\newtheorem{cor}[thm]{Corollary}
\newtheorem{rem}[thm]{Remark}
\newtheorem*{defn}{Definition}

\newcommand{\DD}{\mathbb{D}}
\newcommand{\NN}{\mathbb{N}}
\newcommand{\ZZ}{\mathbb{Z}}
\newcommand{\QQ}{\mathbb{Q}}
\newcommand{\RR}{\mathbb{R}}
\newcommand{\CC}{\mathbb{C}}
\renewcommand{\SS}{\mathbb{S}}

\newcommand{\supp}{\mathop{\mathrm{supp}}}    
\newcommand{\osc}{\mathop{\mathrm{osc}}}    

\newcommand{\re}{\mathop{\mathrm{Re}}}   
\newcommand{\im}{\mathop{\mathrm{Im}}}   
\newcommand{\dist}{\mathop{\mathrm{dist}}}  
\newcommand{\link}{\mathop{\circ\kern-.35em -}}
\newcommand{\spn}{\mathop{\mathrm{span}}}   
\newcommand{\ind}{\mathop{\mathrm{ind}}}   
\newcommand{\rank}{\mathop{\mathrm{rank}}}   
\newcommand{\Fix}{\mathop{\mathrm{Fix}}}   
\newcommand{\codim}{\mathop{\mathrm{codim}}}   
\newcommand{\conv}{\mathop{\mathrm{conv}}}   
\newcommand{\epsi}{\mbox{$\varepsilon$}}
\newcommand{\eps}{\mathchoice{\epsi}{\epsi}
{\mbox{\scriptsize\epsi}}{\mbox{\tiny\epsi}}}
\newcommand{\cl}{\overline}
\newcommand{\pa}{\partial}
\newcommand{\ve}{\varepsilon}
\newcommand{\zi}{\zeta}
\newcommand{\Si}{\Sigma}
\newcommand{\cA}{{\mathcal A}}
\newcommand{\cG}{{\mathcal G}}
\newcommand{\cH}{{\mathcal H}}
\newcommand{\cI}{{\mathcal I}}
\newcommand{\cJ}{{\mathcal J}}
\newcommand{\cK}{{\mathcal K}}
\newcommand{\cL}{{\mathcal L}}
\newcommand{\cN}{{\mathcal N}}
\newcommand{\cR}{{\mathcal R}}
\newcommand{\cS}{{\mathcal S}}
\newcommand{\cT}{{\mathcal T}}
\newcommand{\cU}{{\mathcal U}}
\newcommand{\OM}{\Omega}
\newcommand{\B}{\bullet}
\newcommand{\ol}{\overline}
\newcommand{\ul}{\underline}
\newcommand{\vp}{\varphi}
\newcommand{\AC}{\mathop{\mathrm{AC}}}   
\newcommand{\Lip}{\mathop{\mathrm{Lip}}}   
\newcommand{\es}{\mathop{\mathrm{esssup}}}   
\newcommand{\les}{\mathop{\mathrm{les}}}   
\newcommand{\nid}{\noindent}
\newcommand{\pzr}{\phi^0_R}
\newcommand{\pir}{\phi^\infty_R}
\newcommand{\psr}{\phi^*_R}
\newcommand{\pow}{\frac{N}{N-1}}
\newcommand{\ncl}{\mathop{\mathrm{nc-lim}}}   
\newcommand{\nvl}{\mathop{\mathrm{nv-lim}}}  
\newcommand{\la}{\lambda}
\newcommand{\La}{\Lambda}    
\newcommand{\de}{\delta}    
\newcommand{\fhi}{\varphi} 
\newcommand{\ga}{\gamma}    
\newcommand{\ka}{\kappa}   

\newcommand{\core}{\heartsuit}
\newcommand{\diam}{\mathrm{diam}}

\newcommand{\lan}{\langle}
\newcommand{\ran}{\rangle}
\newcommand{\tr}{\mathop{\mathrm{tr}}}
\newcommand{\diag}{\mathop{\mathrm{diag}}}
\newcommand{\dv}{\mathop{\mathrm{div}}}

\newcommand{\al}{\alpha}
\newcommand{\be}{\beta}
\newcommand{\Om}{\Omega}
\newcommand{\na}{\nabla}

\newcommand{\cC}{\mathcal{C}}
\newcommand{\cM}{\mathcal{M}}
\newcommand{\nr}{\Vert}
\newcommand{\De}{\Delta}
\newcommand{\cX}{\mathcal{X}}
\newcommand{\cP}{\mathcal{P}}
\newcommand{\om}{\omega}
\newcommand{\si}{\sigma}
\newcommand{\te}{\theta}
\newcommand{\Ga}{\Gamma}

\newcommand{\ds}{\displaystyle}

\title[Alexandrov, Serrin, Weinberger, Reilly]{Alexandrov, Serrin, Weinberger, Reilly: \\ simmetry and stability by integral identities}

\author{Rolando Magnanini} 
\address{Dipartimento di Matematica ed Informatica ``U.~Dini'',
Universit\` a di Firenze, viale Morgagni 67/A, 50134 Firenze, Italy.}
    \email{magnanin@math.unifi.it}
    \urladdr{http://web.math.unifi.it/users/magnanin}


\begin{abstract}
The distinguished names in the title have to do with different proofs of the celebrated Soap Bubble Theorem and of radial symmetry in certain overdetermined boundary value problems. We shall give an overeview of those results and indicate some of their ramifications. We will also show how more recent proofs uncover the path to some stability results for the relevant problems.
\end{abstract}

\keywords{Serrin's overdetermined problem, Alexandrov Soap Bubble Theorem, torsional rigidity, constant mean curvature, integral identities, quadrature identities, stability, quantitative estimates}
\subjclass{Primary 35N25, 53A10, 35B35; Secondary 35A23}

\maketitle

\raggedbottom

\section{Introduction}
In this short survey, the author wishes to give an overview of some results related to Alexandrov's Soap Bubble Theorem and Serrin's symmetry result for overdetermined boundary value problems. The presentation will be as untechnical as possible: we shall give no rigorous proofs --- but indicate the relevant references to them --- preferring to focus on ideas and their mutual connections. As the title hints, we will mainly concentrate on the method of integral identities.
\par
We will start by presenting the various proofs of the two results, then we shall explain how they benefit from one another, and hence examine their relatinship to other areas in mathematical analysis. We will finally report on some recent stability results that detail quantitatively how close to the spherical configuration the solution is, if the relevant data are perturbed.

\section{Alexandrov's Soap Bubble Theorem and reflection principle}

Alexandrov's Soap Bubble theorem dates back to $1958$ and states:
\begin{thm}[Soap Bubble Theorem, \cite{Al1}, \cite{Al2}]
\label{th:SBT}
A compact hypersurface, embedded in $\mathbb{R}^N$,  that has constant mean curvature must be a sphere.
\end{thm}
The {\it mean curvature} $H$ of a hypersurface $\cS$ of class $C^2$ at a given point on $\cS$ is the arithmetic mean of its {\it principal curvatures} at that point (see \cite{Re}).
\par
To prove Theorem \ref{th:SBT}, A. D. Alexandrov introduced what is now known as {\it Alexandrov's reflection principle} (see \cite{Al1},\cite{Al2}).
\par
The underlying idea behind Alexandrov's proof is simple: a compact hypersurface $\cS$ is a sphere if and only if it is mirror-symmetric in any fixed direction, that is, for any direction $\te\in\mathbb{S}^{N-1}$, there is a hyperplane $\pi_\te$ orthogonal to $\te$ such that $\cS$ is symmetric in $\pi_\te$.
The technical tools to carry out that idea pertain to the theory of elliptic partial differential equations. To understand why, we give a sketch of Alexandrov's elegant proof.
\par
Let the mean curvature $H$ be constant and suppose by contradiction that $\cS$ is not symmetric in the direction $\te$ (by a rotation, we can always suppose that $\te$ is the upward vertical direction). Then there exists a hyperplane $\pi_\te$ such that at least one of the following occurrences come about: 
\begin{enumerate}[(i)]
\item
the reflection $\cS'$ in $\pi_\te$ of the portion of $\cS$ that stays below $\pi_\te$, touches $\cS$ internally at some point $p\in\cS\setminus\pi_\te$;
\item
$\pi_\te$ is orthogonal to $\cS$ at some point $p\in\cS\cap\pi_\te$.
\end{enumerate}
In both cases, around $p$ we can locally write $\cS$ and $\cS'$ as graphs of two real-valued functions $u$ and $u'$ of $N-1$ variables. If (i) holds, $u$ and $u'$ can be defined on an $(N-1)$-dimensional ball centered at $p$; if (ii) holds instead, $u$ and $u'$ can be defined on an $(N-1)$-dimensional {\it half-ball} centered at $p$, and $p$ belongs to the flat portion of its boundary.  In any case, we have that $u'\le u$ since $\cS'$ stays below $\cS$, and also $u(p)=u'(p)$ and $\na u'(p)=\na u(p)$. 
\par
A partial differential equation now comes about, since both $u$ and $u'$ satisfy the elliptic equation
$$
\frac1{N-1}\,\dv\left(\frac{\na v}{\sqrt{1+|\na v|^2}}\right)=H,
$$
with $H$ constant, being the left-hand side a formula for the mean curvature of the graph of $v$.
A contradiction then occurs because the solutions of that equation satisfy the {\it strong comparison principle} and the {\it Hopf's comparison lemma}. In fact, in case (i), by the strong comparison principle, it should be $u'<u$, whereas we know that $u'(p)=u(p)$; in case (ii), $p$ is on the boundary of the half-ball (the flat part) and hence, by the Hopf's lemma,  it should be that $u_\te'(p)>u_\te(p)$, being $\te$ the normal to the flat part of the boundary of the half-ball. That gives the desired contradiction, since we know that $\na u'(p)=\na u(p)$.
\par
The reflection principle is quite flexible, since its application can be extended to other geometrical settings, such as that of {\it Weingarten's surfaces}, considered by Alexandrov himself.

\section{Serrin's symmetry result and the method of moving planes}
Serrin's symmetry result has to do with certain overdetermined problems for elliptic or parabolic partial differential equations.
In its simplest formulation, it concerns a function $u\in C^1(\ol{\Om})\cap C^2(\Om)$ satisfying the constraints:
\begin{eqnarray}
\label{serrin1}
&\De u=N \ \mbox{ in } \ \Om, \quad u=0 \ \mbox{ on } \ \Ga, \\
\label{serrin2}
&u_\nu=R \ \mbox{ on } \ \Ga.
\end{eqnarray}
Here, $\Om\subset\RR^N$, $N\ge 2$, is a bounded domain with sufficiently smooth, say $C^2$, boundary $\Ga$, $u_\nu$ is the outward normal derivative of $u$ on $\Ga$, and $R$ is a positive constant. 
\par
Since the Dirichlet problem \eqref{serrin1} already admits a {\it unique} solution, the additional requirement \eqref{serrin2} makes the problem {\it overdetermined} and \eqref{serrin1}-\eqref{serrin2} may not admit a solution in general. Thus, the remaining data of the problem --- the domain $\Om$ --- cannot be given arbitrarily.
\par
In fact, Serrin's celebrated symmetry result
states:
\begin{thm}[Radial symmetry, \cite{Se}]
\label{th:Serrin}
The problem \eqref{serrin1}-\eqref{serrin2} admits a solution $u\in C^1(\ol{\Om})\cap C^2(\Om)$
if and only if, up to translations, $\Om$ is a ball of radius $R$ and $u(x)=(|x|^2-R^2)/2$.
\end{thm}
\par
This result inaugurated a new and fruitful field in mathematical research at the confluence of Analysis and Geometry, that has many applications to other areas of mathematics and natural sciences. To be sure, that same result was actually motivated by two concrete problems in Mathematical Physics regarding the torsion of a straight solid bar and the tangential stress of a fluid on the walls of a rectilinear pipe.
\par
The proof given by Serrin in \cite{Se} extends and refines the idea of Alexandrov. In Serrin's setting, the good news is that the hypersurface $\cS\subset\RR^{N+1}$ to be considered is already the graph of a function on $\Om$; however, the bad news is that $\cS$ has now a {\it non-empty boundary}. Moreover, the expected spherical symmetry concerns the base-domain $\Om$, rather then the hypersurface $\cS$: as a matter of fact, Serrin's statement claims that $\cS$ has to be (a portion of) a (spherical) {\it paraboloid} and not a sphere.
\par
In his 1971's proof, J. Serrin brilliantly adapted the reflection principle, by only considering the reflecting  hyperplains $\pi_\te$ orthogonal to {\it horizontal} directions $\te$. The critical occurrences (i) and (ii) then 
take place in a rather modified fashion: in both cases the point $p$ belongs to $\pa\cS=\ol{\cS}\cap(\pa\Om\times\RR)$ and is not a relatively internal point in $\cS$. Thus, the strong comparison principle is ruled out. Nevertheless, in case (i), the Hopf's lemma can still be applied, giving the desired contradiction.
\par
Even so, there is an additional difficulty that one has to deal with in case (ii): the Hopf's comparison lemma can no longer be applied. This is due to the fact that $p$ is not only in $\pa\cS$, but its projection $\ol{p}$ onto $\ol{\Om}$ is placed at a corner on the boundary of the projection $\Om'$ of $\cS'$ onto $\ol{\Om}$ --- $\Om'$ being the domain of the possible application of Hopf's lemma.
\par
To circumvent this obstacle, Serrin established what is now known as {\it Serrin's corner lemma} and concerns the first and second derivatives at $\ol{p}$ of $u'$ and $u$  in the directions $\ell$ entering $\Om'$ from $\ol{p}$: it must hold that either $u'_\ell(\ol{p})<u_\ell(\ol{p})$ or $u'_{\ell\ell}(\ol{p})<u_{\ell\ell}(\ol{p})$ for some $\ell$. After further calculations, this lemma provides the desired contradiction.
\par
This modification of Alexandrov's reflection principle is what is now called the {\it method of moving planes}. The method is very general since, as pointed out by Serrin himself, it applies at least to elliptic equations of the form
$$
a(u,|\na u|)\,\De u+h(u,|\na u|)\,\lan\na^2 u\na u,\na u\ran=f(u,|\na u|),
$$
provided some sufficient conditions are satisfied by the coefficients $a, h$, and $f$ (see \cite{Se} for details) and, more importantly, under the assumption that {\it non-positive solutions} are considered (solutions of \eqref{serrin1} are authomatically negative by the strong maximum principle). Further extensions have also been given during the years by many authors.

\section{Weinberger's proof of Serrin's result}
\label{sec:weinberger}
In the same issue of the journal in which \cite{Se} is published, H.~F.~Weinberger \cite{We} gave a different proof of Theorem \ref{th:Serrin}, based on integration by parts and the Cauchy-Schwarz inequality. 
\par
Weinberger's proof profits of the fact that the so-called {\it P-function} associated to \eqref{serrin1}, defined by
\begin{equation}
\label{def-P}
P=\frac12\,|\na u|^2-u,
\end{equation}
is sub-harmonic in $\Om$, since
\begin{equation}
\label{delta-P}
\De P=|\na^2 u|^2-\frac1{N}\, (\De u)^2\ge 0,
\end{equation}
by the Cauchy-Schwarz inequality applied, for instance, to the two $N^2$-dimensional vectors formed, respectively, by the entries of the identity matrix and those of the {\it hessian matrix} $\na^2 u$. Since $P=R^2/2$ on $\Ga$, then either $P\equiv R^2/2$ or $P<R^2/2$ on $\Om$, by the strong maximum principle. However, the latter occurrence is ruled out by directly calculating that
\begin{equation}
\label{P-Weinberger}
\int_\Om (R^2/2-P)\,dx=0.
\end{equation}
This formula follows by applying the divergence theorem and integration by parts formulas in various forms. Indeed, from \eqref{serrin1} and \eqref{serrin2}, we have:
\begin{equation}
\label{3-identities}
\begin{array}{c}
\ds\int_\Om P\,dx=\left(\frac12+\frac{1}{N}\right)\,\int_\Om |\na u|^2 dx; \\
\ds R\,|\Ga|=\int_\Ga u_\nu\,dS_x=N\,|\Om|;  \quad
(N+2) \int_{\Om} |\na u|^2 \, dx =\int_\Ga u_\nu^2 \,(x\cdot\nu)\,dS_x.
\end{array}
\end{equation}
The second formula sets the correct value for $R$; the third one is a consequence of the well-known  {\it Rellich-Pohozaev identity} (\cite{Po}).
\par
Thus, it must hold that $P\equiv R^2/2$ on $\ol{\Om}$, which implies that $\De P\equiv 0$ in $\Om$.
This means, in turn, that the inequality in \eqref{delta-P} holds with the sign of equality, that is the hessian matrix $\na^2 u$ is proportional to the identity matrix. Then, it is easy to show that  $u$ must equal a {\it quadratic polynomial} of the form
$$
q(x)=\frac12\,(|x-z|^2-a),
$$
for some $z\in\RR^N$ and $a\in\RR$. Since $u=0$ on $\Ga$, we can compute that $a=R^2$, and this implies that $\Ga$ is a sphere centered at $z$ and with radius $R$. 
\par
Weinberger's argument is very elegant, but so far is known to work only for the simple setting \eqref{serrin1}-\eqref{serrin2} or some restricted extensions of it. In particular, it is not known to work if we replace $N$ in \eqref{serrin1} by a non-constant function of $u$.

\section{Reilly's proof of the Soap Bubble Theorem}
In 1982, R.~C.~Reilly found a proof of the Soap Bubble Theorem that bears a resemblance to Weinberger's argument. \par
The key idea is to regard the hypersurface $\cS$ as the zero-level surface of the solution $u$ of \eqref{serrin1}, that is $\cS=\Ga$, and to observe that 
\begin{equation}
\label{Reilly}
\De u =u_{\nu\nu}+(N-1)\,H\,u_\nu \ \mbox{ on } \ \Ga.
\end{equation}
In this formula, that holds on any {\it regular} level surface of a function $u\in C^2(\ol{\Om})$, we agree to still denote by $\nu$ the vector field $\na u/|\na u|$ (that indeed coincides on $\Ga$ with the unit normal field).
\par
As in Section \ref{sec:weinberger}, the radial symmetry of $\Om$ is obtained by showing that $\De P\equiv 0$ in $\Om$. In fact, if $H\equiv H_0$ on $\Ga$ for some constant $H_0$, one can show that
$$
\int_\Om \De P\,dx \le 0.
$$
This inequality follows by using the divergence theorem and by applying a set of formulas similar to \eqref{3-identities}:
\begin{equation*}
\label{3-identities-SBT}
\begin{array}{c}
\ds\int_\Ga P_\nu\,dS_x
=N\,|\Om|-\int_\Ga H\,u_\nu^2\,dS_x; \\
\ds N\,|\Om| H_0=|\Ga|;  \quad
\bigl(N\,|\Om|\bigr)^2\le |\Ga|\,\int_\Ga u_\nu^2\,dS_x.
\end{array}
\end{equation*}
In the first identity, we use \eqref{Reilly} and the divergence theorem. The second formula sets the correct value for the constant $H_0$ and is a consequence of {\it Minkowski's identity},
$$
\int_\Ga H\,(x\cdot\nu)\,dS_x=|\Ga|,
$$
a well-known result in differential geometry (see \cite{Re} for an elementary proof). The last inequality is clearly an application of H\"older's inequality.
\par
Notice that Reilly's argument leaves open the possibility to extend the Soap Bubble Theorem to more general regularity settings, provided a weaker definition of mean curvature is at hand. 
\par
Extensions of Reilly's ideas to the case of the {\it symmetric invariants} of the principal curvatures of $\Ga$ can be found in \cite{Ro}, where a proof of {\it Heintze-Karcher inequality},
$$
\int_\Ga\frac{dS_x}{H}\ge N|\Om|,
$$
is also given, in the same spirit.

\section{The isoperimetric inequality for the torsional rigidity}
\label{sec:torsion}
As an interlude, we present a connection of Serrin's problem \eqref{serrin1}-\eqref{serrin2} to a classical result in shape optimization.
In fact, as also referred to in \cite{Se}, the solution of \eqref{serrin1} has to do with an important quantity in elasticity: the so-called  {\it torsional rigidity} $\tau(\Om)$ of a bar of cross-section $\Om$ (see \cite[pp. 109-119]{So}) that, with the necessary normalizations, can be defined as the maximum of the quotient
$$
Q(v)=\frac{\left(N\,\int_\Om v\,dx\right)^2}{\int_\Om |\na v|^2\,dx},
$$
among all the non-zero functions $v$ in the Sobolev space $W^{1,2}_0(\Om)$. In fact, it turns out that
\begin{equation}
\label{torsional-rigidity}
\tau(\Om)=Q(u)=\int_\Om |\na u|^2\,dx=-N\,\int_\Om u\,dx.
\end{equation}
\par
The following statement is known as the {\it Saint Venant's Principle}: 
\begin{quote}
\it Among sets having given volume the ball maximizes $\tau(\Om)$.
\end{quote}
One proof of this principle hinges on rearrangement tecniques  (see \cite{PSz}).
\par
Here, we shall give an account of the relationship between \eqref{serrin1}-\eqref{serrin2} and the Saint Venant Principle. In fact, once the existence of a maximizing set $\Om_0$ is established, one can show that the solution of \eqref{serrin1} in $\Om_0$ also satisfies \eqref{serrin2} with $\Ga=\pa\Om_0$.
\par 
One way to see that is to introduce the technique of {\it shape derivative}. That consists in hunting for the optimal domain within a one-parameter family $\{\Om_t\}_{t\in\RR}$ of domains that evolve according to a prescribed rule. Thus, if we agree that $\Om_0$ is the domain that maximizes $\tau(\Om_t)$ among all the domains in the family that have prescribed volume $|\Om_t|=V$, then the {\it method of Lagrange multipliers} informs us that there is a number $\la$ such that
$$
T(t)-\la\,[V-V(t)]\le T(0) \ \mbox{ for any } t\in\RR,
$$
and hence 
\begin{equation}
\label{Lagrange}
T'(0)+\la\,V'(0)=0,
\end{equation}
where we mean that $T(t)=\tau(\Om_t)$ and $V(t)=|\Om_t|$.
\par
A convenient way to construct the evolution of the domains $\Om_t$ of the family is to let each of them be the image $\cM_t(\Om)$ of a fixed domain $\Om=\Om_0$ by
a mapping $\cM_t:\RR^N\to\RR^N$ belonging to a family such that:
\begin{equation}
\label{mapping}
\cM_0(x)=x, \quad \cM'_0(x)=\cR(x),
\end{equation}
where the ``prime'' means differentiation with respect to $t$.
\par
Thus, we can consider the solution $u=u(t,x)$ of \eqref{serrin1} in $\Om=\Om_t$ and obtain:
$$
T(t)=-N\,\int_{\Om_t} u(t,x)\,dx \quad \mbox{ and } \quad V(t)=\int_{\Om_t} dx.
$$
The derivatives of $T$ and $V$ can be computed by applying the theory of shape derivatives, stemmed from {\it Hadamard's variational formula} (see \cite[Chapter 5]{HP}). In fact, by a theorem of J.~Liouville, we have that
$$
T'(0)=-N\,\int_{\Om_0} u'(x)\,dx-N\,\int_{\Ga_0} u(x)\,\cR(x)\cdot\nu(x)\,dS_x
$$
and
\begin{equation}
\label{der-V}
V'(0)=\int_{\Ga_0} \cR(x)\cdot\nu(x)\,dS_x,
\end{equation}
where $\Ga_0=\pa\Om_0$ and we have set $u(x)=u(0,x)$ and denoted by $u'(x)$  the derivative of $u(t,x)$ with respect to $t$, evaluated at $t=0$. Moreover, $u'$ turns out to be the solution of the Dirichlet problem
$$
\De u'=0 \ \mbox{ in } \ \Om_0, \quad u'=\na u\cdot\cR \ \mbox{ on } \ \Ga_0.
$$
Also, since $u=0$ on $\Ga_0$ and $u'$ is harmonic in $\Om_0$, we calculate that
\begin{equation}
\label{der-T}
T'(0)=-N\,\int_{\Om_0} u'\,dx=-\int_{\Om_0} u'\,\De u\,dx=-\int_{\Ga_0} u'\,u_\nu\,dS_x,
\end{equation}
after an application of Gauss-Green's formula.
\par
Next, we choose
$$
\cR(x)=\phi(x)\,\nu(x),
$$
where $\phi$ is any compactly supported continuous function and $\nu$ is a proper extension of the unit normal vector field to a tubular neighborhood of $\Ga$ (for instance the choice $\nu(x)=\na\de_\Ga(x)$, where $\de_\Ga(x)$ is the distance of $x$ from $\Ga$, will do).
\par
Therefore, by this choice of $\cR$, putting together \eqref{Lagrange}, \eqref{der-V}, and \eqref{der-T} gives that
$$
\int_{\Ga_0} (u_\nu^2-\la)\,\phi\,dS_x=0.
$$
Since $\phi$ is arbitrary, we infer that $u_\nu^2\equiv\la$ on $\Ga_0$ and we compute $\la=R^2$.
\par
Theorem \ref{th:Serrin} thus confirms Saint Venant's principle. A sufficient regularity assumption that guarantees that this argument runs is that $\Ga_0$ is locally the graph of a differentiable function with H\"older continuous derivatives.

\section{Dual formulation and quadrature domains}

In \cite{PS}, the following characterization is proved.

\begin{thm}[Dual formulation: harmonic domain]
\label{th:Payne-Schaefer}
A function $u\in C^1(\ol{\Om})\cap C^2(\Om)$ is solution of \eqref{serrin1}-\eqref{serrin2} if and only if
the following mean value property 
\begin{equation}
\label{dual-formulation}
\frac1{|\Om|}\int_\Om h\,dx=\frac1{|\Ga|}\int_\Ga h\,dS_x 
\end{equation}
holds for any harmonic function $h\in C^0(\ol{\Om})\cap C^2(\Om)$.
\end{thm}
The proof is a straightforward consequence of Gauss-Green's formula for laplacians. 
\par
A domain $\Om$ 
such that \eqref{dual-formulation} holds for any harmonic function is named a {\it harmonic domain} (see \cite{RS}). Thus, the following corollary ensues.

\begin{cor}
\label{cor:PS}
The euclidean ball is the only bounded harmonic domain in $\RR^N$.
\end{cor}

Theorem \ref{th:Payne-Schaefer} and Corollary \ref{cor:PS} are due to L.~Payne and P.~W.~Schaefer \cite{PS}, who also provide a proof of Theorem \ref{th:Serrin} that modifies Weinberger's argument and gets rid of the use of the maximum principle for $P$.

In \cite{PS}, \eqref{dual-formulation} is regarded as a dual formulation of the overdetermined problem \eqref{serrin1}-\eqref{serrin2}, since it entails a linear functional,
$$
\cH(\Om)\ni h\mapsto L(h)=\frac1{|\Om|}\int_\Om h\,dx-\frac1{|\Ga|}\int_\Ga h\,dS_x,
$$
defined on the set $\cH(\Om)$ of functions in $C^0(\ol{\Om})\cap C^2(\Om)$ that are harmonic in $\Om$.
\par
The identity \eqref{dual-formulation} recalls the well-known Gauss {\it mean value theorems} for harmonic functions: if $\Om$ is a ball and $p$ its center, then
$$
h(p)=\frac1{|\Om|}\int_\Om h\,dx \ \mbox{ and } \ h(p)=\frac1{|\Ga|}\int_\Ga h\,dS_x 
$$
for any harmonic function $h\in C^0(\ol{\Om})\cap C^2(\Om)$. It is interesting to remark that each mean value property characterizes the ball (or the sphere) with respect to the class of harmonic functions, as \eqref{dual-formulation} does (see \cite{PS}, \cite{Ku}).

\medskip

The dual formulation can also be connected to the theory of {\it quadrature domains} introduced by D.~Aharonov and B.~Gustafsson in the 1970's (see \cite{AS}, \cite{GS}). A bounded domain $\Om$ in the complex plane $\CC$ is a (classical) quadrature domain if there exist finitely many points $p_1,\dots, p_m\in\Om$ and coefficients $c_{jk}\in\CC$ so that
\begin{equation}
\label{quadrature}
\int_\Om f(z)\,dx dy=\sum_{j=1}^m \sum_{k=0}^{n_j} c_{jk} f^{(k)}(p_j)
\end{equation}
for any {\it holomorphic function} $f$ in $\Om$; here, $f^{(0)}=f$ and $f^{(k)}$ denotes the $k$-th derivative of $f$.
\par
Formula \eqref{quadrature} is called a {\it quadrature identity}.  For instance, the ball centered at $p$ is a quadrature domain that corresponds to $m=1$, $n_1=0$ and $p_1=p$.
The term ``quadrature'' thus refers to the fact that, in a quadrature domain $\Om$, \eqref{quadrature} provides an exact {\it quadrature formula} to compute the integral on the left-hand side.
A remarkable fact is that quadrature domains have applications to problems in Mathematical Physics such as the Hele-Shaw problem in fluid dynamics and other free-boundary and/or  inverse problems (see \cite{GS} and the references therin). 
\par
The notion of quadrature domain can be extended in two ways. One can replace the finite sum in \eqref{quadrature} by the integral
$$
\int f(z)\,d\mu,
$$
where $\mu$ is some {\it signed measure} (e.g., in \eqref{quadrature} $\mu$ would be the linear combination of Dirac deltas and their derivatives at given points).
Moreover, one can establish a generalization to higher-dimensional domains, by replacing holomorphic functions by harmonic functions --- in this case, the domain is often called a {\it harmonic quadrature domain}.
Therefore, a harmonic domain --- that is satisfying \eqref{dual-formulation} for any $h\in\cH(\Om)$ --- is a harmonic quadrature domain relative to (a multiple of) the {\it surface measure} on $\Ga$.

\section{Fundamental identities}

In this section, we will show that Weinberger's and Reilly's proofs can be further refined to encode all the information given in Serrin's and Alexandrov's problems into two identities. In the following two results, we use the defininitions:
\begin{equation}
\label{def-R-H0}
R=\frac{N |\Om|}{|\Ga|}, \quad H_0=\frac{|\Ga|}{N |\Om|} \ \mbox{ and } \ q(x)=\frac12\,|x-z|^2-a, \ x\in\RR^N,
\end{equation}
where $z\in\RR^N$ and $a\in\RR$ are given parameters.

\begin{thm}[Fundamental identity for Serrin's problem]
\label{th:wps}
Let $\Om \subset \mathbb R^N$ be a bounded domain with boundary $\Ga$ of class $C^{1,\al}$, $0<\al\le 1$.
\par
Then, the solution $u$ of  \eqref{serrin1} satisfies identity:
\begin{equation}
\label{idwps}
\int_{\Om} (-u)\,\left\{ |\na ^2 u|^2- \frac{ (\De u)^2}{N} \right\}\,dx=
\frac{1}{2}\,\int_\Ga \left( u_\nu^2 - R^2 \right) \,(u_\nu-q_\nu)\,dS_x.
\end{equation}
\end{thm}

The identity is announced in \cite{MP1} and proved in \cite{MP2}. Its proof is obtained by polishing the arguments in \cite{We} and \cite{PS} and juggling around with integration by parts. 

\begin{thm}[Fundamental identity for the Soap Bubble Theorem]
Let $\Om$ be a bounded domain with boundary $\Ga$ of class $C^2$ and let $u$ be the solution of \eqref{serrin1}.
\par
Then, it holds true that
\begin{multline}
\label{identity-SBT}
\frac1{N-1}\int_{\Om} \left\{ |\na ^2 u|^2-\frac{(\De u)^2}{N}\right\}dx+
\frac1{R}\,\int_\Ga (u_\nu-R)^2 dS_x = \\
\int_{\Ga}(H_0-H)\,(u_\nu-q_\nu)\,u_\nu\,dS_x+
\int_{\Ga}(H_0-H)\, (u_\nu-R)\,q_\nu\, dS_x.
\end{multline}
\end{thm}

The identity is proved in \cite{MP2} by slightly modifying one that was proved in \cite{MP1}. Its proof is obtained by polishing the argument in \cite{Re}.
\par
From identity \eqref{idwps} it is clear that, if the right-hand side is zero --- and that surely occurs if \eqref{serrin2} is in force --- then $\De P\equiv 0$ owing to \eqref{delta-P}, and hence radial symmetry ensues, as already observed. 
\par
Since both summands at the left-hand side are non-negative, the same conclusion results from \eqref{identity-SBT}, if its right-hand side is null --- and that holds if $H$ is constant. It should also be noticed that $H\equiv H_0$ implies independently that $u$ satisfies \eqref{serrin2}.

\medskip

One more comment is in order. If we turn back to Section \ref{sec:torsion}, we see that
$$
T'(0)+R^2\,V'(0)=\frac12\,(u_\nu^2-R^2)\,(u_\nu-q_\nu),
$$
if we choose $\phi=(q_\nu-u_\nu)/2$, and hence \eqref{idwps} can be written as:
\begin{equation}
\label{idwps-2}
\int_{\Om} (-u)\,\left\{ |\na ^2 u|^2- \frac{ (\De u)^2}{N} \right\}\,dx=
\int_\Ga [T'(0)+R^2\,V'(0)]\,dS_x.
\end{equation}
\par
Ergo, this identity tells us something more about the Saint Venant Principle.
\begin{thm}
\label{th:priviledged-flow}
A domain $\Om$ is a ball if the function
$$
\RR\ni t\mapsto \tau(\Om_t)+R^2\,(V-|\Om_t|)
$$
obtained by modifing $\Om$ by the rule \eqref{mapping} with $\cR=\cR^*$ and
$$
\cR^*=\frac12\,(q_\nu-u_\nu)\,\nu
$$
has a critical point at $t=0$.
\end{thm}
Actually, it is enough that the derivative is non-positive at $t=0$. Thus, the flow generated by that $\cR^*$ is quite a priviledged one. 
\par
Theorem \ref{th:priviledged-flow} seems to be new.

\section{Stability: in the wake of Alexandrov and Serrin}

In this and the next section, we will present recent results on the stability for the radial configuration in the Soap Bubble Theorem and Serrin's problem. Roughly speaking, the question is how much a hypersurface $\Ga$ is near a sphere, if its mean curvature $H$ --- or, alternatively, the normal derivative on $\Ga$ of the solution $u$ of \eqref{serrin1} --- is near a constant in some norm.
Technically speaking, one may look for two concentric balls $B_{\rho_i}$ and $B_{\rho_e}$, with radii $\rho_i$ and $\rho_e$, $\rho_i<\rho_e$, such that
\begin{equation}
\label{balls}
\Ga\subset \ol{B}_{\rho_e}\setminus B_{\rho_i} 
\end{equation} 
 and
\begin{equation}
\label{stability}
\rho_e-\rho_i\le \psi(\eta),
\end{equation}
where
$\psi:[0,\infty)\to[0,\infty)$ is a continuous function vanishing at $0$ and $\eta$ is a suitable measure of the deviation of $u_\nu$ or $H$ from being a constant.

\par
In this section, I will briefly give an account of the results in this direction obtained by means of quantitative versions of the reflection principle or the method of moving planes. 
\par
The problem of stability for Serrin's problem has been considered for the first time in \cite{ABR}. There, for a $C^{2,\al}$-regular domain $\Om$, it is proved that, if $u$ is the solution of \eqref{serrin1}, there exist constants $C, \ve> 0$ such that \eqref{balls} and \eqref{stability} hold for 
$$
\psi(\eta)=C\,|\log \eta|^{-1/N} \quad \mbox{ and } \quad \eta=\nr u_\nu-c\nr_{C^1(\Ga)}<\ve,
$$
for some constant $c$. It is important to observe that the validity of that inequality extend to the case in which at the right-hand side of the Poisson's equation in \eqref{serrin1} the number $N$ is replaced by a locally Lipschitz continuous function $f(u)$. In this case, only positive solutions are considered and the constants $C, \ve$
also depend on $f$ and the regularity of $\Ga$.
\par
In the same general framework, the stability estimate of \cite{ABR} has been improved in \cite{CMV}. There, it is in fact shown that \eqref{balls} and \eqref{stability} hold for 
$$
\psi(\eta)=C\,\eta^\tau \quad \mbox{ and } \quad 
\eta=\sup_{\substack{x,y \in \Ga\\ \ x \neq y}} \frac{|u_\nu(x) - u_\nu(y)|}{|x-y|}<\ve.
$$
The exponent $\tau\in (0,1)$ can be computed for a general setting and, if $\Om$ is convex,
is proved to be arbitrarily close to $1/(N+1)$. 
\par
The only quantitative estimate for symmetry in the Soap Bubble Theorem, based on Alexandrov's reflection principle, is proved in \cite{CV} and is optimal. There, it is shown that, if $\Ga$ is an $N$-dimensional, $C^2$-regular, connected, closed hypersurface embedded in $\RR^{N}$, there exist constants $C, \ve> 0$ such that \eqref{balls} and \eqref{stability} hold for
$$
\psi(\eta)=C\,\eta \quad \mbox{ and } \quad \eta=\max_\Ga H-\min_\Ga H<\ve.
$$ 
The two constants depend on $N$, upper bounds for the principal curvatures of $\Ga$, and $|\Ga|$.
The result is optimal, because is attained for ellipsoids.
\medskip

We conclude this section by giving a brief outline of the arguments used to obtain stability for Serrin's problem. The arguments are substantially those of \cite{ABR}, that have been refined in \cite{CMS} and \cite{CMV}, and adapted to the situation of the Soap Bubble Theorem in \cite{CV}.
\par
The idea is to fix a direction $\te$ and define an approximate set $X_\eta\subset\Om$, mirror-symmetric with respect to a hyperplane $\pi_\te$ orthogonal to $\te$, that fits $\Om$ ``well'', in the sense that it is the maximal $\te$-symmetric set contained in $\Om$ and such that its order of approximation of $\Om$ can be controlled by $C\,\psi(\eta)$. It turns out that this approximation process does not depend on the particular direction $\te$ chosen. Thus, one defines an approximate center of symmetry $p$ as the intersection of $N$ mutually orthogonal hyperplanes $\pi_{\te_1}, \dots, \pi_{\te_N}$ of symmetry. It then becomes apparent that, in any other direction $\te$, the approximation deteriorates only by replacing $C$ by a possibly larger constant. It is thus possible to define the desired balls in \eqref{balls} by centering them at $p$.
\par
Technically, the control of the approximation by $C\,\psi(\eta)$ is made possible by the application of {\it Harnack's inequality} and {\it Carleson's (or boundary Harnack's) inequality}, which are the quantitative versions of the already mentioned maximum principle, Hopf's lemma, and Serrin's corner lemma. The improvement obtained in \cite{CMV} is the result of a refinement of Harnack's inequality in suitable cones.

\section{Stability: in the wake of Reilly and Weinberger}

Quantitative inequalities for the Soap Bubble Theorem and Serrin's problem can also be obtained by following the tracks of Reilly's and Weinberger's proofs of symmetry. 
\par
In \cite{CM}, based on the proof of Heintze-Karcher's inequality given in \cite{Ro}, it is shown that \eqref{balls}-\eqref{stability} hold for
$$
\psi(\eta)=C\,\eta^\frac1{2(N+1)} \quad \mbox{ and } \quad \eta=\max_\Ga |H_0-H|<\ve, 
$$
for some positive constants $C, \ve$; $\ve$ should be sufficiently small so as to guarantee that $\Ga$ is strictly {\it mean convex} (that means that $H>0$ on $\Ga$) --- the realm of validity of Heintze-Karcher's inequality. The exponent is not optimal; however, an estimate is also given in \cite{CM} that gives a finer description of hypersurfaces having their mean curvature close to a constant. In fact, such an estimate specifies how $\Ga$ can be close to the boundary of a disjoint union of balls.
\par
That result has been improved in various directions in \cite{MP1}. In fact, based on a version of identity \eqref{identity-SBT}, \eqref{balls}-\eqref{stability} are shown to hold for some positive constants $C, \ve$ and
$$
\psi(\eta)=C\,\eta^{\tau_N} \quad \mbox{ and } \quad \eta=\int_\Ga (H_0-H)^+\,dS_x<\ve, 
$$
where $\tau_N=1/2$ for $N=2, 3$ and $\tau_N=1/(N+2)$ for $N\ge 4$. Here we mean $(t)^+=\max(t,0)$ for $t\in\RR$. 
\par
That approximation is not restricted to the class of strictly mean convex hypersurfaces, but is valid for $C^2$-regular hypersurfaces, as in \cite{CV}. Differently from \cite{CV} and \cite{CM}, it replaces the uniform deviation from $H_0$ by a weaker average deviation and yet, compared to \cite{CM}, it improves the relevant stability exponent. 
\par
A further advance has been recently obtained in \cite{MP2}. In fact, it holds that
\begin{equation}
\label{optimal-stability-SBT}
\rho_e-\rho_i\le C\,\nr H_0-H\nr_{2,\Ga}^{\tau_N} \quad \mbox{ if } \quad \nr H_0-H\nr_{2,\Ga}<\ve,
\end{equation}
where $\tau_N=1$ for $N=2, 3$ and $\tau_N=2/(N+2)$ for $N\ge 4$ --- that is the stability exponent doubles.
Therefore, for $N=2, 3$ an optimal Lipschitz inequality (as in \cite{CV}) is established, {\it even for a weaker average deviation}; it seems realistic to expect that a similar Lipschitz estimate holds also for $N\ge 4$.
\par
In \cite{MP2}, an inequality involving the following slight modification of the so-called {\it Fraenkel asymmetry},
\begin{equation*}
\label{asymmetry}
\cA(\Om)=\inf\left\{\frac{|\Om\De B^x|}{|B^x|}: x \mbox{ center of a ball $B^x$ with radius $R$} \right\},
\end{equation*}
has also been proved:
$$
\cA(\Om) \le C\,\nr H_0-H\nr_{2,\Ga}.
$$
Here, $\Om\De B^x$ denotes the symmetric difference of $\Om$ and $B^x$, and $R$ is the constant defined in \eqref{def-R-H0}. This inequality holds for any $N\ge 2$. Under sufficient assumptions, the number $\cA(\Om)$ can be linked to the difference $\rho_e-\rho_i$ (see \cite{MP2}); however, the resulting stability inequality is poorer than \eqref{optimal-stability-SBT}.
\par
At the end of this section, we shall explain how these results have been made possible by parallel estimates for Serrin's problem. 

\medskip

The first improvement of the logarithmic estimate obtained in \cite{ABR} for problem \eqref{serrin1}-\eqref{serrin2} has been given in \cite{BNST}. There, the idea of working on integral identities and inequalities has also been put in action for the first time. By a combination of the ideas of Weinberger and the use of (pointwise) {\it Newton's inequalities} for the hessian matrix $\na^2 u$ of $u$, the solution $u$ of \eqref{serrin1} is shown to satisfy \eqref{balls}-\eqref{stability} for
$$
\psi(\eta)=C\,\eta^{\tau_N} \quad \mbox{ and } \quad \eta=\max_\Ga |u_\nu-c|<\ve,
$$
where $c$ is some reference constant.  In \cite{BNST}, it is also considered the possibility to measure the deviation of $u_\nu$ from a constant by the $L^1$-norm, that is with $\eta=\nr u_\nu-c\nr_{1,\Ga}$ and,
by assuming an appropriate a priori bound for $|\na u|$ on $\Ga$, it is shown that $\Om$ can be approximated in measure by a finite number of mutually disjoint balls $B_i$. The obtained exponent is $\tau_N=1/(4N+9)$.
\par
Recently, that approach has been greatly improved in \cite{Fe} where, rather than by \eqref{balls} and \eqref{stability}, the closeness of $\Om$ to a ball is measured by the asymmetry $\cA(\Om)$:
$$
\cA(\Om)\le C\,\nr u_\nu-R\nr_{2,\Ga}.
$$
\par
Turning back to an approximation of type \eqref{balls}-\eqref{stability}, the following formula is derived in \cite{MP2}:
\begin{equation}
\label{stability-L2-Serrin}
\rho_e-\rho_i\le C\, \nr u_\nu-R\nr_{2,\Ga}^\frac2{N+2} \quad \mbox{ if } \quad \nr u_\nu-R\nr_{2,\Ga}<\ve.
\end{equation}

This inequality clearly makes better than those in \cite{BNST} and \cite{CMV}, even if we replace the uniform norm or the Lipschitz semi-norm by an $L^2$-deviation.

\medskip

As promised, we conclude this section by giving an outline of the proof of \eqref{stability-L2-Serrin} and \eqref{optimal-stability-SBT},
by drawing the reader's attention on how the results benefit from the sharp formulas \eqref{idwps} and \eqref{identity-SBT} and their interaction.
\par
To simplify matters, it is convenient to re-write \eqref{idwps} in terms of the {\it harmonic} function $h=q-u$: it holds that
\begin{equation}
\label{idwps-h}
\int_{\Om} (-u)\, |\na ^2 h|^2\,dx=
\frac{1}{2}\,\int_\Ga ( R^2-u_\nu^2)\, h_\nu\,dS_x.
\end{equation}
Notice that $h=q$ on $\Ga$ and hence the {\it oscillation} of $h$ on $\Ga$ can be bounded from below by $\rho_e-\rho_i$:
\begin{equation*}
\label{oscillation}
\max_{\Ga} h-\min_{\Ga} h=\frac12\,(\rho_e^2-\rho_i^2)\ge \frac12\,(|\Om|/|B|)^{1/N}(\rho_e-\rho_i).
\end{equation*}
Thus, \eqref{balls}-\eqref{stability} will be obtained if can bound that oscillation in terms of the left-hand side of \eqref{idwps-h}. In fact, its right-hand side can be easily bounded in terms of the desired $L^2$-deviation of $u_\nu$ from $R$.
\par
To carry out this plan, the following inequalities, proved in \cite[Lemma 3.3]{MP1}, are decisive:
$$
\max_{\Ga} h-\min_{\Ga} h \le C\,\left(\int_\Om h^2 dx\right)^\frac1{N+2}\le C\,\left(\int_\Om |\na h|^2 dx\right)^\frac1{N+2}.
$$
In the first inequality, it is crucial that $h$ is harmonic in $\Om$, because the mean value property for harmonic functions on balls is used; the second inequality follows from an application of the Poincar\'e inequality, since we can choose $a$ so that $h$ has average zero on $\Om$. \par
Next, notice that the obtained inequalities for the oscillation of $h$ and \eqref{idwps-h} do not depend on the particular choice of $z\in\RR^N$. A good choice for $z$ turns out to be a minimum (or any critical) point of $u$, so that it is guaranteed that $z\in\Om$. That choice has yet a more important benefit: since now $\na h(z)=0$, the Hardy-Poincar\`e-type inequality
\begin{equation}
\label{boas-straube}
\int_\Om v(x)^2 dx\le C \int_\Om (-u) |\na v(x)|^2 dx,
\end{equation}
that holds for any harmonic function $v$ that is zero at some point, can be applied to each first partial derivative of $h$ so as to eventually obtain that 
$$
\max_{\Ga} h-\min_{\Ga} h \le C\,\left(\int_\Om (-u) |\na^2 h|^2 dx\right)^\frac1{N+2}.
$$
Inequality \eqref{boas-straube} can be derived by using an inequality proved in \cite{BS} or \cite{HS} (see \cite{Fe} or \cite{MP2} for details).
\par
The last inequality and \eqref{idwps-h} easily give a stability bound in terms of the deviation $\eta=\nr u_\nu-R\nr_{1,\Ga}$. Nevertheless, one can gain a better estimate by observing that,
if $u_\nu-R$ tends to $0$, also $h_\nu$ does. Quantitavely, this fact can be expressed by the inequality
\begin{equation}
\label{Feldman}
\nr h_\nu\nr_{2,\Ga}\le C\,\nr u_\nu-R\nr_{2,\Ga},
\end{equation}
that can be derived from \cite{Fe}. 
Thus, \eqref{stability-L2-Serrin} will follow by using this inequality, after an application of H\"older's inequality to the right-hand side of \eqref{idwps-h}.
\par
By keeping track of the constants $C$ and $\ve$ in the various inequalities, one can show that those in 
\eqref{stability-L2-Serrin} only depend on $N$, the diameter of $\Om$ and the radii of the optimal interior and exterior touching balls to $\Ga$ (see \cite{MP2}).

\medskip

In order to prove \eqref{optimal-stability-SBT}, we use the new identity \eqref{identity-SBT}, that also reads as:
\begin{multline}
\label{identity-SBT-h}
\frac1{N-1}\int_{\Om} |\na ^2 h|^2 dx+
\frac1{R}\,\int_\Ga (u_\nu-R)^2 dS_x = \\
-\int_{\Ga}(H_0-H)\,h_\nu\,u_\nu\,dS_x+
\int_{\Ga}(H_0-H)\, (u_\nu-R)\,q_\nu\, dS_x.
\end{multline}
In fact, discarding the first summand at its left-hand side and applying H\"older's inequality to the two terms at its right-hand side and \eqref{Feldman} thereafter, yields that
$$
\nr u_\nu-R\nr_{2,\Ga}\le C\,\nr H_0-H\nr_{2,\Ga}.
$$
Inequality \eqref{optimal-stability-SBT} then follows again from \eqref{identity-SBT-h} and the estimate 
$$
\max_{\Ga} h-\min_{\Ga} h\le C \left(\int_{\Om} |\na ^2 h|^2\,dx\right)^{\tau_N/2}
$$
already obtained in \cite{MP1}. 

\section{Stability for a harmonic domain}

We conclude this paper by deriving a stability inequality for the mean value property established in \eqref{dual-formulation}, in the spirit of that given in \cite{CFL} for the classical Gauss mean value property for harmonic functions on balls in $\RR^N$.
\par
We start by writing an identity,
$$
\frac1{|\Om|}\int_\Om h\,dx-\frac1{|\Ga|}\int_\Ga h\,dS_x=\frac1{N |\Om|}\,\int_\Ga h\,(u_\nu-R)\,dS_x,
$$
that holds at least for any function $h\in C^2(\Om)\cap C^0(\ol{\Om})$ which is harmonic in $\Om$, if $u$ is the solution of \eqref{serrin1}. The identity can be easily obtained by an application of Gauss-Green's formula 
and the use of \eqref{serrin1}.
\par
The measure of the deviation of the two mean values from one another can be obtained by taking the norm of the linear functional $L$ defined at the left-hand side of the identity on some relevant normed space. We may for instance consider the {\it Hardy-type space} $\cH^p(\Om)$ of harmonic functions in $\Om$, whose trace on $\Ga$ is a function in $L^p(\Ga)$. Thus, we compute that
\begin{multline*}
\nr L\nr_p=\sup\left\{\left|\frac1{|\Om|}\int_\Om h\,dx-\frac1{|\Ga|}\int_\Ga h\,dS_x\right|: h\in\cH^p(\Om), \nr h\nr_{p,\Ga}\le 1\right\}=\\
\frac1{N |\Om|}\,\nr u_\nu-R\nr_{p',\Ga},
\end{multline*}
where $p'$, as usual, is the conjugate exponent of $p$. 
These computations work for any $p\in[1,\infty]$. 
\par
Therefore, by choosing $p=2$, we obtain the inequality
$$
\rho_e-\rho_i\le C\,\nr L\nr_2^\frac2{N+2} \quad \mbox{ if } \quad \nr L\nr_2<\ve
$$
from \eqref{stability-L2-Serrin}. This inequality is new.

\end{document}